\documentclass[12pt]{amsart}
\usepackage{amssymb,latexsym}

\numberwithin{equation}{section}

\newtheorem{theorem}{Theorem}[section]

\newtheorem{conjecture}[theorem]{Conjecture}
\newtheorem{corollary}[theorem]{Corollary}

\theoremstyle{definition}

\newtheorem{remark}[theorem]{Remark}
\newtheorem{example}[theorem]{Example}
\newtheorem{definition}[theorem]{Definition}
\newtheorem{problem}[theorem]{Problem}

\def\AA{\mathcal{A}}
\def\ZZ{\mathbb{Z}}
\def\CC{\mathbb{C}}

\def\RR{\mathbb{R}}
\def\QQ{\mathbb{Q}}
\def\cc{\mathbf{c}}
\def\ii{\mathbf{i}}
\def\xx{\mathbf{x}}

\def\FFcal{\mathcal{F}}

\def\ex{\mathbf{ex}}
\def\Upper{\mathcal{U}}

\renewcommand{\eqref}[1]{{\rm (\ref{#1})}}
\newcommand{\mat}[4]{\left(\!\!\begin{array}{cc}
#1 & #2 \\ #3 & #4 \\
\end{array}\!\!\right)}

\begin{document}

\title[Cluster algebras]
{Cluster algebras: notes for 2004 IMCC \\[.1in]
(Chonju, Korea, August 2004)}

\author{Andrei Zelevinsky}
\address{\noindent Department of Mathematics, Northeastern University,
 Boston, MA 02115, USA}
\email{andrei@neu.edu}

\begin{abstract}
This is an expanded version of the notes for the two lectures
at the 2004 International Mathematics Conference
(Chonbuk National University, August 4--6, 2004).
The first lecture discusses the origins of cluster algebras, with
the focus on total positivity and geometry of double Bruhat cells
in semisimple groups.
The second lecture introduces cluster algebras and discusses
some basic results, open questions and conjectures.
\end{abstract}

\date{July 23, 2004}

 \thanks{Research supported in part
by NSF grant DMS-0200299.}



\maketitle

\tableofcontents

\section{Introduction}

Cluster algebras, introduced and studied in~\cite{ca1,ca2,ca3},
are a class of axiomatically defined commutative
rings equipped with a distinguished set of generators (cluster
variables) grouped into overlapping subsets (clusters) of the same
finite cardinality.
The original motivation for this theory was to design an
algebraic framework for total positivity
and canonical bases in semisimple algebraic groups.
Since its inception, the theory of cluster algebras has developed
several interesting connections and applications:
\begin{itemize}
\item
Discrete dynamical systems based on rational
recurrences~\cite{fz-laurent, carroll-speyer, speyer}.
\item
$Y$-systems in thermodynamic
Bethe Ansatz~\cite{yga}.
\item
A new family of convex polytopes (generalized associahedra)
associated to finite root systems~\cite{yga,cfz}.
%including as special cases Stasheff's associahedron, and
%Bott-Taubes cyclohedron.
\item
Quiver representations~\cite{bmrrt, bmr, ccs, mrz}.
\item
Grassmannians, projective configurations and
their tropical analogues~\cite{scott,speyer-williams}.
\item
Quantum cluster algebras, Poisson geometry and Teichm\"uller theory~\cite
{fock-goncharov1, fock-goncharov2, gsv1, gsv2, qca}.
\end{itemize}

In these lectures we concentrate on the following aspects of the
theory of cluster algebras.
Lecture~1 discusses the origins of the theory,
with the focus on total positivity and geometry of double Bruhat cells
in semisimple groups (a closely related connection to canonical bases
was discussed in~\cite{camblec}).
Lecture~2 introduces cluster algebras
and discusses some basic results, open problems and conjectures.

I thank Sergey Fomin for helpful comments and suggestions.

\section{Lecture~1: Total positivity and double Bruhat cells}
\label{sec:lecture1}

A matrix is \emph{totally positive} \ (resp.\ \emph{totally nonnegative})
if all its minors are positive (resp.\ nonnegative) real numbers.
The first systematic study of these classes of matrices was undertaken
in the 1930s by I.~J.~Schoenberg, and by F.~R.~Gantmacher and M.~G.~Krein
(the references to their papers and further connections and
applications of totally positive matrices can be found
in~\cite{fz-intel}).
The interest in this subject intensified in the last decade
due to the discovery by G.~Lusztig~\cite{lusztig}
of a surprising connection between total positivity and
canonical bases for quantum groups.
Among other things, he extended the subject by defining
the \emph{totally positive variety} $G_{> 0}$ and the
\emph{totally nonnegative variety} $G_{\geq 0}$
inside every complex reductive group~$G$.
These ideas were further developed
in~\cite{bfz-tp, bz-schubert, fz-double, fz-osc}.
We are going to present some of the results obtained in these papers.
For technical reasons, we restrict ourselves to the case when~$G$
is semisimple and simply connected.
(The results on classical total positivity are recovered when
$G = SL_{r+1}(\CC)$ is the group of complex
$(r+1) \times (r+1)$ matrices with determinant~$1$.)

Lusztig defined $G_{> 0}$ and $G_{\geq 0}$ parametrically.
As shown in~\cite{fz-osc}, Lusztig's original definition is
equivalent to setting
\begin{align}
\label{eq:gtp-gtnn}
G_{> 0}& = \{x \in G: \Delta_{\gamma,\delta}(x) > 0
\,\, \text{for all} \,\, \gamma,\delta\}, \\
\nonumber
G_{\geq 0} &= \{x \in G: \Delta_{\gamma,\delta}(x) \geq 0
\,\, \text{for all} \,\, \gamma,\delta\},
\end{align}
for a family of regular functions~$\Delta_{\gamma,\delta}$
on~$G$ called \emph{generalized minors}
(they specialize to ordinary minors in the case $G = SL_{r+1}(\CC)$).
To be more specific, let $\omega_1, \dots, \omega_r$ be the
fundamental weights of~$G$, and~$W$ the Weyl group
(the standard terminology and notation used here is explained in more detail
in~\cite{fz-double}).
The indices~$\gamma$ and~$\delta$ of a generalized minor are two
elements of the $W$-orbit of the same~$\omega_i$, i.e., two
extremal weights in the same fundamental representation~$V_{\omega_i}$ of~$G$;
the corresponding minor $\Delta_{\gamma,\delta}$ is a suitably
normalized matrix element of~$V_{\omega_i}$ associated with
weights~$\gamma$ and~$\delta$.
(For $G = SL_{r+1}(\CC)$, the Weyl group $W$ is identified
with the symmetric group~$S_{r+1}$, and
$V_{\omega_i} = \bigwedge^i \CC^{r+1}$, the $i$-th exterior power
of the standard representation.
All the weights of $V_{\omega_i}$
are extremal, and they are in bijection with the $i$-subsets of
$[1,r+1] = \{1, 2, \dots, r+1\}$, so that~$W = S_{r+1}$
acts on them in a natural way, and $\omega_i$ corresponds to~$[1,i]$.
If~$\gamma$ and~$\delta$ correspond to $i$-subsets $I$ and $J$, respectively,
then $\Delta_{\gamma,\delta} = \Delta_{I,J}$ is the minor with
the row set~$I$ and the column set~$J$.)

The natural geometric framework for studying the varieties
$G_{> 0}$ and $G_{\geq 0}$ is provided by
\emph{double Bruhat cells}
\begin{equation}
\label{eq:double-cell}
G^{u,v} = B u B \cap B_- v B_- \,;
\end{equation}
here $u,v \in W$, and~$B$ and~$B_-$ are two opposite Borel
subgroups in~$G$.
(For $G = SL_{r+1}(\CC)$, the standard choice for~$B$
(resp.~$B_-$) is the subgroup of upper (resp.~lower) triangular
matrices.)
The group~$G$ has two \emph{Bruhat decompositions}, with respect
to~$B$ and $B_-\,$:
\begin{equation}
G = \bigcup_{u \in W} B u B = \bigcup_{v \in W} B_- v B_- \ .
\end{equation}
Therefore, $G$ is the disjoint union of all double Bruhat cells.
The term ``cell" is somewhat misleading since the topology of
$G^{u,v}$ may be non-trivial.
According to \cite[Theorem~1.1]{fz-double}, the variety
$G^{u,v}$ is isomorphic to a Zariski open subset of an affine space of dimension
$r+\ell(u)+\ell(v)$, where~$r$ is the rank of~$G$,
and~$w \mapsto \ell(w)$ is the usual length function on~$W$.
(Recall that~$W$ is a Coxeter group generated by \emph{simple
reflections} $s_1, \dots, s_r$.
For $w \in W$, the length $\ell = \ell(w)$ is the shortest length of a
sequence of indices $\ii = (i_1, \dots, i_\ell)$ such that
$w = s_{i_1} \cdots s_{i_\ell}$; such a sequence is called a
\emph{reduced word} for~$w$.
For $G = SL_{r+1}(\CC)$, the simple reflection
$s_i \in W = S_{r+1}$ is the
transposition of~$i$ and~$i+1$,
and $\ell(w)$ is the number of inversions of a permutation~$w$.)

Following \cite{fz-double}, we define the \emph{totally positive part}
of $G^{u,v}$ by setting
\begin{equation}
G^{u,v}_{> 0} = G^{u,v} \cap G_{\geq 0} \ .
\end{equation}
For instance, if $w_\circ$ is the longest element of~$W$
then $G^{w_\circ, w_\circ}$ is the open double Bruhat cell given by:
\begin{equation}
G^{w_\circ, w_\circ} = \{x \in G:
\Delta_{\omega_i, w_\circ \omega_i}(x) \neq 0, \,\,
\Delta_{w_\circ \omega_i, \omega_i}(x) \neq 0 \,\, \text{for all}
\,\, i \in [1,r]\};
\end{equation}
and we have
\begin{equation}
\label{eq:tp-part-biggest-cell}
G^{w_\circ, w_\circ}_{> 0} = G_{> 0} \, .
\end{equation}
(As a nice exercise in linear algebra, the reader is invited to
check \eqref{eq:tp-part-biggest-cell} for $G = SL_3(\CC)$; in
other words, if a matrix $x \in SL_3(\CC)$ has all its minors
nonnegative, and the minors
$$\Delta_{1,3}(x) = x_{13}, \,\, \Delta_{12,23}(x), \,\,
\Delta_{3,1}(x) = x_{31}, \,\, \Delta_{23,12}(x)$$
are non-zero, then all the minors of~$x$ are positive.)

The structure of~$G^{u,v}_{> 0}$ was described by
G.~Lusztig~\cite{lusztig}: it
turns out to be isomorphic to $\RR_{> 0}^{r+\ell(u)+\ell(v)}$.
To establish this, Lusztig introduced local coordinates in each
$G^{u,v}$ which consist of rational functions that are not necessarily regular.
To sharpen this result, we introduce the following definition
which is a slight variation of the one given in \cite{fz-double}.

\begin{definition}
\label{def:tp-basis}
A \emph{TP-basis} for $G^{u,v}$ is a a collection of regular
functions $F= \{f_1, \dots, f_m\} \subset \CC[G^{u,v}]$ with the
following properties:
\begin{enumerate}
  \item The functions $f_i$ are algebraically independent and
  generate the field of rational functions $\CC(G^{u,v})$; in
  particular, $|F| = m = r+\ell(u)+\ell(v)$.
  \item The mapping $(f_1, \dots, f_m): G^{u,v} \to \CC^m$ restricts to a
  biregular isomorphism $U(F) \to \CC_{\neq 0}^m$,
  where $U(F)$ is the locus of all $x \in G^{u,v}$ such that
$f_k(x) \neq 0$ for all $k \in [1,m]$.
  \item The mapping $(f_1, \dots, f_m): G^{u,v} \to \CC^m$
  restricts to an isomorphism $G^{u,v}_{> 0} \to \RR_{> 0}^m$.
\end{enumerate}
\end{definition}

The following theorem from~\cite{fz-double} extends the
results previously obtained in~\cite{bfz-tp, bz-schubert}.

\begin{theorem}
\label{th:tp-basis}
Every reduced word~$\ii$ for $(u,v) \in W \times W$ gives rise to a
certain TP-basis $F_\ii$ for~$G^{u,v}$ (an explicit description
of $F_\ii$ is given below).
\end{theorem}

The set $F_\ii$ is a collection of generalized minors given explicitly as follows.
For technical reasons that will soon become clear, we will represent~$\ii$
as a sequence of indices $(i_1, \dots, i_m)$
from the set $-[1,r] \cup [1,r]$, where $m = r + \ell(u) + \ell(v)$,
the set $-[1,r]$ consists of indices $\{-1, \dots, -r\}$, we have
$i_j = j$ for $j \in [1,r]$, and
\begin{equation}
\label{eq:ii-products}
s_{-i_{r+1}} \cdots s_{-i_m} = u, \quad
s_{i_{r+1}} \cdots s_{i_m} = v \, ,
\end{equation}
with the convention
\begin{equation}
\label{eq:s-negative}
s_{-i} = 1 \,\, \text{for} \,\, i \in [1,r] \, .
\end{equation}
In this notation, we have
\begin{equation}
\label{eq:Fi}
F_\ii =
%\{\} \, \cup \,
\{\Delta_{\gamma_k, \delta_k}: k \in [1,m]\},
\end{equation}
where
\begin{equation}
\label{eq:gamma-k}
\gamma_k = s_{-i_1} \cdots s_{-i_k} \omega_{|i_k|}, \quad
\delta_k = s_{i_m} \cdots s_{i_{k+1}} \omega_{|i_k|} \, .
\end{equation}

\begin{example}
\label{example:SL3-w0w0}
Let $G=SL_3(\CC)$, and let $u=v=w_\circ=s_1s_2s_1=s_2s_1s_2$ be the order-reversing
permutation (the element of maximal length in the symmetric group~$W=S_3$).
Take $\ii=(1,2,1,2,1,-1,-2,-1)$.
Then we have
\begin{equation}
\label{eq:f-sl3}
F_\ii = \{\Delta_{1,3}, \Delta_{12,23}, \Delta_{1,2},
\Delta_{12,12}, \Delta_{1,1}, \Delta_{2,1}, \Delta_{23,12},
\Delta_{3,1}\}.
\end{equation}
\end{example}

Theorem~\ref{th:tp-basis} can be seen as the first step towards introducing
the cluster algebra structure in the coordinate ring
$\CC[G^{u,v}]$ (the collection $F_\ii$ will give rise to a
cluster in the future theory).
An impetus for introducing this structure came from realizing that, even when
a reduced word~$\ii$ is allowed to vary, the corresponding
families~$F_\ii$ may not be sufficient for answering some natural
questions about double Bruhat cells.
Here is one such question: describe and enumerate the connected
components of the real part of~$G^{u,v}$.
An important special case of this problem when $u = 1$ and
$v = w_\circ$ for $G=SL_{r+1}(\CC)$ was solved in~\cite{ssv1,ssv2}
(see also related work~\cite{rie1,rie2}).
The general case was handled in \cite{z-imrn} using the results
and ideas from \cite{ssvz} (for follow-up see~\cite{gsv0,seven}).
The main idea of the solution in \cite{z-imrn} was to replace
$G^{u,v}$ with a ``simpler" Zariski open subvariety~$U$ such that
the codimension in $G^{u,v}$ of the complement of~$U$ is greater
than~$1$; if furthermore, $U$ is ``compatible" with the real part
of $G^{u,v}$, then replacing  $G^{u,v}$ by~$U$ does not change
the structure of connected components of the real part.
In many cases (including the one covered in~\cite{ssv1,ssv2}) one
can take as~$U$ the union of all open sets of
the form $U(F_\ii)$, where~$\ii$ runs over reduced words for~$(u,v)$.
However, in general the complement of the latter open subvariety has
codimension~$1$, so one needs something else.

Here is the main construction from~\cite{z-imrn}.
Fix a pair~$(u,v) \in W \times W$, and a reduced word
$\ii = (i_1, \dots, i_m)$ for~$(u,v)$ as above.
For $k \in [1,m]$, abbreviate~$f_k = \Delta_{\gamma_k, \delta_k}$,
where $\gamma_k$ and $\delta_k$ are given by \eqref{eq:gamma-k}.
We say that an index~$k \in [1,m]$ is \emph{$\ii$-exchangeable} if
$r < k \leq m$, and $|i_p| = |i_k|$ for some $p > k$.
Denote by $\ex = \ex_\ii \subset [1,m]$ the subset of
$\ii$-exchangeable indices.
Note that the subset
$$\{f_k: k \in [1,m] - \ex\} \, \subset \, F_\ii$$
depends only on~$u$ and~$v$, not on the particular choice of the reduced word~$\ii$.
(Indeed, $f_j = \Delta_{\omega_j, v^{-1}\omega_j}$
for $j \in [1,r]$; and $f_k = \Delta_{u \omega_i,\omega_i}$
if~$k$ is the last occurrence of $\pm i$ in~$\ii$.)
In particular, the cardinality~$n$ of~$\ex$ depends only on~$u$ and~$v$.
We also have the following important property:
\begin{equation}
\label{eq:coeff-not-vanish}
\text{for $k \in [1,m] - \ex$, the function $f_k$ vanishes nowhere on~$G^{u,v}$.}
\end{equation}

\begin{theorem}
\label{th:adjacent-clusters}
There is an integer $m \times n$ matrix $\tilde B = \tilde B(\ii) = (b_{ik})$
with rows labeled by~$[1,m]$ and columns labeled by~$\ex$,
satisfying the following properties:
\begin{enumerate}
\item
For every~$k \in \ex$, the function
\begin{equation}
\label{eq:exchange-f}
f'_k = \frac
{\prod\limits_{b_{ik}>0} f_i^{b_{ik}}+
\prod\limits_{b_{ik}<0} f_i^{-b_{ik}}}{f_k} \, .
\end{equation}
is regular on~$G^{u,v}$, i.e., belongs to $\CC[G^{u,v}]$.
\item
For every~$k \in \ex$, the collection
\begin{equation}
\label{eq:Fik}
F_{\ii;k} = F_\ii - \{f_k\} \cup \{f'_k\} \, \subset \, \CC[G^{u,v}]
\end{equation}
obtained from~$F_\ii$ by exchanging~$f_k$ with~$f'_k$ is a
TP-basis for $G^{u,v}$.
\item
Let~$U \subset G^{u,v}$ be the Zariski open subset given by
$$U = U(F_\ii) \, \cup \, \bigcup_{k \in \ex} U(F_{\ii;k}).$$
Then the complement of~$U$ in $G^{u,v}$ has codimension greater
than~$1$.
\end{enumerate}
\end{theorem}

The matrix~$\tilde B = \tilde B(\ii)$ in
Theorem~\ref{th:adjacent-clusters} is given explicitly as follows
(cf.~\cite[(8.7)]{qca}): for $p \in [1,m]$ and $k \in \ex$, we have
\begin{equation}
\label{eq:tildeB-entries}
b_{pk} = b_{pk}(\ii) =
\begin{cases}
-\varepsilon (i_k) & \text{if $p = k^-$;} \\
-\varepsilon (i_k) a_{|i_p|, |i_k|} &
\text{if $p < k < p^+ < k^+, \, \varepsilon (i_k) = \varepsilon(i_{p^+})$,}\\
 & \text{or $p < k < k^+ < p^+, \, \varepsilon (i_k) = -\varepsilon (i_{k^+})$;} \\
\varepsilon (i_p) a_{|i_p|, |i_k|} &
\text{if $k < p < k^+ < p^+, \, \varepsilon (i_p) = \varepsilon(i_{k^+})$,}\\
& \text{or $k < p < p^+ < k^+, \, \varepsilon (i_p) = -\varepsilon (i_{p^+})$;} \\
\varepsilon (i_p) & \text{if $p = k^+$;} \\
0 & \text{otherwise,}
\end{cases}
\end{equation}
where we use the following notation and conventions:
\begin{itemize}
\item
$\varepsilon(i) = \pm 1$ for $i \in \pm [1,r]$;
\item
$A = (a_{ij})_{i,j \in [1,r]}$ is the Cartan matrix of~$G$,
that is, the transition matrix between simple roots and
fundamental weights:
\begin{equation}
\label{eq:cartan-transition}
\alpha_j = \sum_{i=1}^r a_{ij} \omega_i \ .
\end{equation}
\item
for $k \in [1,m]$, we denote by $k^+ = k^+_\ii$
the smallest index $\ell$ such that $k < \ell \leq m$ and
$|i_\ell| = |i_k|$; if $|i_k| \neq |i_\ell|$ for $k < \ell \leq m$, then
we set $k^+ = m+1$.
\item
$k^- = k^-_\ii$ denotes the index $\ell$ such that $\ell^+=k$;
if such an $\ell$ does not exist, we set $k^- = 0$.
\end{itemize}

\begin{example}
\label{ex:sl3-exchanges}
Let~$G$ and $\ii$ be the same as in
Example~\ref{example:SL3-w0w0}.
Then~$\tilde B = \tilde B(\ii)$ is an $8 \times 4$ matrix given by
$$\tilde B =
  \begin{pmatrix}
-1 &  0 &  0 &  0\\
1 &  -1 &  0 &  0\\
0 & 1 &  -1 &  0\\
-1 &  0 & 1 &  -1\\
1 &  -1 &  0 & 1\\
0 & 1 &  -1 &  0\\
0 &  -1 &  0 & 1\\
0 &  0 &  0 &  -1
  \end{pmatrix},$$
where the columns are indexed by the set $\ex = \{3,4,5,6\}$
(this matrix given in \cite[Example~3.2]{qca} is a slight modification
of the one in \cite[Example~2.5]{ca3}).
Applying \eqref{eq:exchange-f} to the functions in
\eqref{eq:f-sl3}, we obtain
\begin{align*}
f'_3 &= \frac{f_2 f_5 + f_1 f_4}{f_3} =
\frac{(x_{12} x_{23} - x_{13} x_{22})x_{11} +
x_{13}(x_{11} x_{22} - x_{12} x_{21})}{x_{12}}
= \Delta_{12,23} \, ,\\
f'_4 &= \frac{f_3 f_6 + f_2 f_5 f_7}{f_4} =
x_{12}x_{21}x_{33}-x_{12}x_{23}x_{31}-
x_{13}x_{21}x_{32}+x_{13}x_{22}x_{31} \, ,\\
f'_5 &= \frac{f_4 + f_3 f_6}{f_5} = x_{22} \, ,\\
f'_6 &= \frac{f_5 f_7 + f_4 f_8}{f_6} =
\Delta_{13,12} \, .
\end{align*}
Note that, while each of the families $F_{\ii;3}, F_{\ii;5}$,
and~$F_{\ii;6}$ is of the form $F_{\ii'}$ for some reduced
word~$\ii'$, this is no longer true for $F_{\ii;4}$, if only for
the reason that~$f'_4$ is not even a minor.
%This makes the claim that $F_{\ii;4}$ is a TP-basis more intriguing.
\end{example}

In view of Theorem~\ref{th:adjacent-clusters}~(3),
we have $\CC[G^{u,v}] = \CC[U]$, which
implies the following description of the coordinate ring $\CC[G^{u,v}]$.

\begin{corollary}
\label{cor:CGuv}
The subalgebra $\CC[G^{u,v}]$ of the field of rational functions
$\CC(G^{u,v})$ is the intersection of~$n+1$ Laurent polynomial rings:
\begin{equation}
\label{eq:CGuv-Laurent}
\CC[G^{u,v}] = \CC[F_\ii^{\pm 1}] \cap \bigcap_{k \in \ex}
\CC[F_{\ii;k}^{\pm 1}],
\end{equation}
where the collections~$F_{\ii;k}$ are given by \eqref{eq:Fik}.
\end{corollary}

Motivated by Corollary~\ref{cor:CGuv}, we now introduce
\emph{upper cluster algebras}, which first appeared in~\cite{ca3}.
Our definition will not be the most general: as in~\cite{qca},
we restrict our attention to cluster algebras of \emph{geometric type}.

Let~$m$ and~$n$ be two positive integers with $m \geq n$.
Let~$\FFcal$ be the field of rational
functions over~$\QQ$ in~$m$ independent (commuting) variables.

\begin{definition}
\label{def:seed}
A \emph{seed} (of geometric type) in $\FFcal$ is a pair
$(\tilde \xx, \tilde B)$, where
\begin{itemize}

\item $\tilde \xx = \{x_1, \dots, x_m\}$ is an
algebraically independent subset of~$\FFcal$ which generates~$\FFcal$.

\item $\tilde B$ is an $m \times n$ integer matrix with rows
labeled by $[1,m]$
and columns labeled by an $n$-element subset $\ex \subset [1,m]$,
satisfying the two conditions:
\begin{enumerate}
\item
the $n \times n$ submatrix of $\tilde B$ with rows and columns
labeled by $\ex$ is \emph{skew-symmetrizable}, i.e.,
$d_i b_{ik} = -d_k b_{ki}$ for some positive
integers~$d_i \,\, (i, k \in \ex)$.
\item $\tilde B$ has full rank~$n$.
\end{enumerate}
\end{itemize}
The seeds are defined up to a relabeling of elements of $\tilde \xx$
together with the corresponding relabeling of rows and columns of $\tilde B$.
\end{definition}

As shown in \cite[Proposition~2.6]{ca3}, the matrix~$\tilde B =
\tilde B(\ii)$ given by \eqref{eq:tildeB-entries} satisfies the
conditions in Definition~\ref{def:seed}.

In analogy with \eqref{eq:exchange-f} and
\eqref{eq:Fik}, for every~$k \in \ex$, we set
\begin{equation}
\label{eq:exchange}
x'_k = \frac
{\prod\limits_{b_{ik}>0} x_i^{b_{ik}}+
\prod\limits_{b_{ik}<0} x_i^{-b_{ik}}}{x_k} \in \FFcal\, ,
\end{equation}
and
\begin{equation}
\label{eq:adjacent-cluster}
\tilde \xx_k = \tilde \xx - \{x_k\} \cup \{x'_k\} \, .
\end{equation}

\begin{definition}
\label{def:upper-cluster-algebra}
The \emph{upper cluster algebra}
$\Upper(\tilde \xx, \tilde B) = \Upper(\tilde B)$
is the subring of the ambient field~$\FFcal$ given by
\begin{equation}
\label{eq:upper}
\Upper(\tilde \xx, \tilde B) = \Upper(\tilde B)
 = \ZZ[\tilde \xx^{\pm 1}] \cap \bigcap_{k \in \ex}
\ZZ[\tilde \xx_k^{\pm 1}],
\end{equation}
where $\ZZ[\tilde \xx^{\pm 1}]$ stands for the ring of integer
Laurent polynomials in the variables from~$\tilde \xx$.
\end{definition}

In this terminology, Corollary~\ref{cor:CGuv} takes the following
form (cf.~\cite[Theorem~2.10]{ca3}):

\begin{corollary}
\label{cor:CGuv-upper-cluster}
The coordinate ring $\CC[G^{u,v}]$ of any double Bruhat cell
is naturally isomorphic to the complexification of the upper
cluster algebra $\Upper(\tilde B(\ii))$
associated with the matrix $\tilde B(\ii)$ given by
\eqref{eq:tildeB-entries}.
\end{corollary}

To make this abstraction useful, we need to develop some tools for
better understanding of upper cluster algebras.
This does not seem to be an easy task.
In particular, it is not at all clear how big the intersection
of several Laurent polynomial rings in \eqref{eq:upper} can be.
For instance, the following problem raised in \cite[Problem~1.27]{ca3}
is still open:

\begin{problem}
\label{prob:finite-gen-upper}
For which matrices~$\tilde B$ the algebra~$\Upper(\tilde B)$
is finitely generated?
\end{problem}

In the next lecture we introduce a mechanism
for generating a lot of new elements (\emph{cluster
variables}) in~$\Upper(\tilde B)$.

\section{Lecture~2: Some results and conjectures on cluster algebras}

With the motivation outlined above,
we now move towards introducing cluster algebras.
We retain the terminology and notation in
Definition~\ref{def:upper-cluster-algebra}.
In addition, we denote
$\xx = \{x_j: j \in \ex\} \subset \tilde \xx$, and
$\cc = \tilde \xx - \xx$.
We refer to the indices from $\ex$ as \emph{exchangeable} indices,
to $\xx$ as the \emph{cluster}, and to
the $n \times n$ submatrix $B$ of $\tilde B$ with rows and columns
labeled by~$\ex$ as the \emph{principal part} of~$\tilde B$ and
the \emph{exchange matrix} of a seed
$(\tilde \xx, \tilde B)$.

Following \cite[Definition~4.2]{ca1}, we say that the
$m \times n$ matrix $\tilde B'$ is obtained from $\tilde B$ by
\emph{matrix mutation} in direction~$k \in \ex$,
and write $\tilde B' = \mu_k (\tilde B)$
if the entries of $\tilde B'$ are given by
\begin{equation}
\label{eq:matrix-mutation}
b'_{ij} =
\begin{cases}
-b_{ij} & \text{if $i=k$ or $j=k$;} \\[.05in]
b_{ij} + \displaystyle\frac{|b_{ik}| b_{kj} +
b_{ik} |b_{kj}|}{2} & \text{otherwise.}
\end{cases}
\end{equation}
One can show (see e.g., \cite[Proposition~2.3]{qca})
that matrix mutations preserve the conditions in
Definition~\ref{def:seed}.
Furthermore, the principal part of $\tilde B'$ is equal to $\mu_k(B)$,
and $\mu_k$ is involutive: $\mu_k (\tilde B') = \tilde B$.

\begin{definition}
\label{def:seed-mutation}
Let $(\tilde \xx, \tilde B)$ be a seed in $\FFcal$.
For any exchangeable index $k$, the \emph{seed mutation} in direction~$k$
transforms $(\tilde \xx, \tilde B)$ into the seed
$\mu_k(\tilde \xx, \tilde B)=(\tilde \xx', \tilde B')$, where
$\tilde \xx' = \tilde \xx_k$ is given by
\eqref{eq:exchange}--\eqref{eq:adjacent-cluster}, and $\tilde B' = \mu_k(\tilde B)$.
\end{definition}

Note that $(\tilde \xx', \tilde B')$ is indeed a seed, and
the seed mutation is involutive, i.e.,
$\mu_k (\tilde \xx', \tilde B') = (\tilde \xx, \tilde B)$.
Therefore, we can define an equivalence relation on seeds as follows:
we say that $(\tilde \xx, \tilde B)$ is mutation-equivalent to
$(\tilde \xx', \tilde B')$ and write
$(\tilde \xx, \tilde B) \sim (\tilde \xx', \tilde B')$ if $(\tilde \xx', \tilde B')$
can be obtained from $(\tilde \xx, \tilde B)$ by a sequence of seed mutations.
Note that all seeds $(\tilde \xx', \tilde B')$ mutation-equivalent to
a given seed $(\tilde \xx, \tilde B)$ share the same set
$\cc = \tilde \xx - \xx$.
%Let $\ZZ[\cc^{\pm 1}] \subset \FFcal$ be the ring of integer
%Laurent polynomials in the elements of~$\cc$.
%We call the set $\xx' = \tilde \xx' - \cc$ the \emph{cluster}
%of a seed $(\tilde \xx', \tilde B')$.
Let $\mathcal{X} = \mathcal{X}(\mathcal{S})$ denote
the union of clusters of all seeds in $\mathcal{S}$
(recall that the cluster of $(\tilde \xx', \tilde B')$ is
the subset $\xx' = \tilde \xx' - \cc \subset \tilde \xx'$).
We refer to the elements of $\mathcal{X}$ as \emph{cluster variables}.

Now everything is in place for defining cluster algebras.

\begin{definition}
\label{def:cluster-algebra}
Let $\mathcal{S}$ be a mutation-equivalence class of seeds.
The \emph{cluster algebra} $\AA(\mathcal{S})$ associated with $\mathcal{S}$
is the $\ZZ[\cc^{\pm 1}]$-subalgebra of the ambient field $\FFcal$
generated by all cluster variables: thus, we have
$\AA(\mathcal{S}) = \ZZ[\cc^{\pm 1}, \mathcal{X}]$.
The cardinality~$n$ of every cluster is called the
\emph{rank} of $\AA(\mathcal{S})$.
\end{definition}

Since $\mathcal{S}$ is uniquely determined by each of its seeds
$(\tilde \xx, \tilde B)$, we sometimes
denote $\AA(\mathcal{S})$ as $\AA(\tilde \xx,\tilde B)$, or even simply $\AA(\tilde B)$,
because $\tilde B$ determines this algebra uniquely up to an
automorphism of~$\FFcal$.

%Definitions~\ref{def:seed-mutation} and \ref{def:cluster-algebra} are justified by
The following result first appeared in~\cite[Theorem~3.1]{ca1}.

\begin{theorem}[Laurent phenomenon]
\label{th:A-in-U}
The cluster algebra~$\AA(\mathcal{S})$ is contained in the upper
cluster algebra $\Upper(\tilde \xx,\tilde B)$ for every seed
$(\tilde \xx, \tilde B)$ in~$\mathcal{S}$.
Equivalently,
$\AA(\mathcal{S}) \ \subset \ \ZZ[\tilde \xx^{\pm 1}]$,
i.e., every element of $\AA(\mathcal{S})$ is an integer
Laurent polynomial in the variables from~$\tilde \xx$.
\end{theorem}

A new proof of Theorem~\ref{th:A-in-U} was given in \cite{ca3},
where it was obtained as a consequence of the following result.

\begin{theorem}[\cite{ca3}, Theorem~1.5]
\label{th:upper-mutation-invariant}
The upper cluster algebra $\Upper(\tilde \xx,\tilde B)$ is
mutation-invariant, that is,
$\Upper(\tilde \xx,\tilde B)$ and $\Upper(\tilde \xx',\tilde B')$
coincide whenever the seeds $(\tilde \xx,\tilde B)$ and $(\tilde \xx',\tilde B')$
are mutation-equivalent.
\end{theorem}

In view of Theorem~\ref{th:upper-mutation-invariant}, we can use the notation
$\Upper(\mathcal{S}) = \Upper(\tilde \xx,\tilde B)$, where
$\mathcal{S}$ is the mutation-equivalence class of a seed
$(\tilde \xx,\tilde B)$.

It is worth mentioning that Theorem~\ref{th:upper-mutation-invariant} and the
resulting Laurent phenomenon carry over to the quantum setting
developed in~\cite{qca}.

We now present some open problems and conjectures on cluster algebras.
The following open problem was raised in \cite[Problem~1.25]{ca3}.

\begin{problem}
\label{problem:when-upper=cluster}
When is a cluster algebra equal to the corresponding upper cluster
algebra?
\end{problem}

A nice sufficient condition for the equality
$\AA(\mathcal{S}) = \Upper(\mathcal{S})$ was found in~\cite{ca3}.
Following~\cite{ca3}, we call a seed~$(\tilde \xx,\tilde B)$
%(as well as the corresponding cluster algebra)
\emph{acyclic} if there is a linear ordering of~$\ex$ such that
$b_{ij} \geq 0$ for all $i,j \in \ex$ with $i < j$ (since the
principal part~$B$ of~$\tilde B$ is skew-symmetrizable, we also have $b_{ji} \leq 0$).

\begin{theorem}[\cite{ca3}, Theorem~1.18, Corollary~1.19]
\label{th:acyclic}
If~$\mathcal{S}$ contains  an acyclic seed then $\AA(\mathcal{S}) = \Upper(\mathcal{S})$.
Moreover, if a seed~$(\tilde \xx,\tilde B)$ is acyclic then
$\AA(\tilde \xx,\tilde B) = \Upper(\tilde \xx,\tilde B)
= \ZZ[\cc^{\pm 1}, x_k, x'_k \, (k \in \ex)]$, where
$x'_k$ is given by \eqref{eq:exchange}; in particular, this
algebra is finitely generated.
\end{theorem}

An important combinatorial invariant of a cluster algebra
$\AA(\mathcal{S})$ is its \emph{exchange graph}.
The vertices of this graph correspond to the seeds
from a mutation-equivalence class~$\mathcal{S}$,
and the edges correspond to seed mutations.
By definition, the exchange graph is connected and $n$-regular
(that is, every vertex has degree~$n$).
Besides these properties, not much is known about the exchange
graph in general.
Here are a few conjectures about this graph (some of them already
appeared in \cite{ca2}).

\begin{conjecture}
\label{con:exchange-general}
\begin{enumerate}
\item
The exchange graph of an algebra $\AA(\tilde \xx,\tilde B)$
depends only on the principal part~$B$ of the matrix~$\tilde B$.
\item
Every seed in~$\mathcal{S}$ is uniquely determined by its cluster;
thus, the exchange graph can be thought of as having the clusters
as vertices, with two clusters adjacent if and only if their
intersection has cardinality~$n-1$.
\item
For every cluster variable $x \in \mathcal{X}$, the seeds
from~$\mathcal{S}$ whose clusters contain~$x$ form a
\underline{connected} subgraph of the exchange graph.
\item
The seeds from~$\mathcal{S}$ whose exchange matrix~$B$ is acyclic
form a \underline{connected} subgraph (possibly empty) of the exchange graph.
\end{enumerate}
\end{conjecture}

Another important problem is to find a nice combinatorial labeling
of cluster variables.
More generally, we would like to parametrize cluster monomials
defined as follows.

\begin{definition}
\label{def:cluster-monomial}
A \emph{cluster monomial} is a monomial in the cluster variables
such that all the variables occurring in it belong to the same cluster.
\end{definition}

As a labeling set for cluster monomials, we would like to use the
denominators in the Laurent expansion with respect to a given cluster.
To be more specific, we fix a cluster~$\xx$ of some seed from~$\mathcal{S}$,
and choose a numbering of its elements: $\xx = \{x_1, \dots, x_n\}$.
Every non-zero element $y \in \Upper(\mathcal{S})$ can be uniquely
written as
\begin{equation}
\label{eq:Laurent-normal-form}
y = \frac{P(x_1, \dots, x_n)}{x_1^{d_1} \cdots x_n^{d_n}} \, ,
\end{equation}
where $P(x_1, \dots, x_n)$ is a polynomial with coefficients in~$\ZZ[\cc^{\pm 1}]$
which is not divisible by any cluster variable from~$\xx$.
We denote
\begin{equation}
\label{eq:denominator-vector}
\delta(y) = \delta_\xx(y) = (d_1, \dots, d_n) \in \ZZ^n,
\end{equation}
and call the integer vector~$\delta(y)$ the \emph{denominator vector}
of~$y$ with respect to the cluster~$\xx$.
For instance, the elements of~$\xx$ and their immediate exchange
partners have the denominator vectors
\begin{equation}
\label{eq:denom-xj}
\delta (x_j) = - e_j, \quad \delta (x'_j) = e_j \quad (j \in [1,n]),
\end{equation}
where $e_1, \dots, e_n$ are the standard basis vectors in~$\ZZ^n$.
Note also that the map $y \mapsto \delta(y)$ satisfies:
\begin{equation}
\label{eq:delta-multiplicative}
\delta(yz) = \delta(y) + \delta(z) \ .
\end{equation}

\begin{conjecture}
\label{con:denominators}
Distinct cluster monomials have distinct denominator vectors with
respect to any cluster.
%, implying that the cluster monomials are
%linearly independent over $\ZZ[\cc^{\pm 1}]$.
\end{conjecture}

Our last group of questions/conjectures addresses \emph{positivity}.
We define the \emph{positive cones} $\Upper(\mathcal{S})_{\geq 0}$
and $\AA(\mathcal{S})_{\geq 0}$ by setting
\begin{equation}
\label{eq:positive-cone}
\Upper(\mathcal{S})_{\geq 0} = \bigcap_{(\tilde \xx,\tilde B) \in \mathcal{S}}
\ZZ_{\geq 0}[\tilde \xx^{\pm 1}], \quad
\AA(\mathcal{S})_{\geq 0} =
\AA(\mathcal{S}) \ \cap \ \Upper(\mathcal{S})_{\geq 0}\, ;
\end{equation}
that is, $\Upper(\mathcal{S})_{\geq 0}$ (resp. $\AA(\mathcal{S})_{\geq 0}$)
consists of all elements $y \in \Upper(\mathcal{S})$
(resp. $y \in \AA(\mathcal{S})$) such that, for every
seed~$(\tilde \xx,\tilde B) \in \mathcal{S}$, the Laurent expansion
of~$y$ in the variables from~$\tilde \xx$ has nonnegative
coefficients.
A non-zero element of a positive cone will be referred to as \emph{positive}.

\begin{problem}
\label{prob:positive-cones}
Describe the positive cones $\Upper(\mathcal{S})_{\geq 0}$
and $\AA(\mathcal{S})_{\geq 0}$.
\end{problem}

To be more specific: let us call a positive element (of either of the
positive cones) \emph{indecomposable} if it cannot be written as
a sum of two positive elements.

\begin{problem}
\label{prob:indec-positive}
Describe the set of all indecomposable elements
in either of the cones $\Upper(\mathcal{S})_{\geq 0}$
and $\AA(\mathcal{S})_{\geq 0}$.
\end{problem}

\begin{conjecture}
\label{con:cluster-positive}
\begin{enumerate}
\item
(\cite{ca1})
Every cluster variable is a positive element of the cluster
algebra~$\AA$, i.e., its Laurent expansion
in the variables from~$\tilde \xx$ has nonnegative
coefficients for every seed $(\tilde \xx,\tilde B)$ of~$\AA$ .
\item
Every cluster monomial is an indecomposable positive element.
\end{enumerate}
\end{conjecture}

The above conjectures are non-trivial already in rank~$2$.
Conjectures~\ref{con:exchange-general} and \ref{con:denominators}
hold for any cluster algebra of rank~$2$ (the former follows from
the results in \cite{ca1}, while the latter is proved in~\cite{sz}).
Conjecture~\ref{con:cluster-positive} is still open even in
rank~$2$; in~\cite{sz}, it was proved for special rank~$2$
cluster algebras with the exchange matrix of the form
\begin{equation}
\label{eq:B-rank2}
B = \mat{0}{b}{-c}{0}  \quad (b,c > 0, \,\, bc \leq 4) \, .
\end{equation}
Furthermore, it was shown in \cite{sz} that in the case of \eqref{eq:B-rank2},
the indecomposable positive elements form a $\ZZ$-basis of the cluster algebra.
This raises the following general question.

\begin{problem}
\label{prob:indec-positive-basis}
For which matrices~$\tilde B$,
the indecomposable positive elements form a $\ZZ$-basis
of the corresponding cluster algebra?
\end{problem}

All the above problems and conjectures become much more tractable
for the finite type cluster algebras defined as follows.

\begin{definition}
\label{def:finite-type}
A cluster algebra $\AA(\mathcal{S})$ is of
\emph{finite type} if the mutation-equivalence class~$\mathcal{S}$
is finite, i.e., the exchange graph has finitely many vertices.
\end{definition}

A classification of cluster algebras of finite type was given in \cite{ca2}.
Remarkably, this classification turns out to be identical to the famous Cartan-Killing classification of
semisimple Lie algebras and finite root systems.
To be more specific: for an $n\times n$ Cartan matrix $A\!=\!(a_{ij})$ of finite type
(cf.~\eqref{eq:cartan-transition}),
define a skew-symmetrizable  matrix $B(A)\!=\!(b_{ij})$ by
\[
b_{ij}=
\begin{cases}
\ \ 0 & \text{if $i=j$;} \\
\varepsilon(i)\, a_{ij} & \text{if $i\neq j$,}
\end{cases}
\]
where the sign function $\varepsilon: [1,n] \to \{1,-1\}$ is such that
$a_{ij}< 0 \ \Longrightarrow\ \varepsilon(i)=-\varepsilon(j)$
(the existence of such a sign function follows from the fact that
every Dynkin diagram is a tree, hence a bipartite graph).

\begin{theorem}[\cite{ca2}]
\label{th:fin-type-class}
A cluster algebra $\AA$ is of finite type
if and only if the exchange matrix at some
seed of~$\AA$ is of the form~$B(A)$, where $A$ is a Cartan matrix
of finite type.
Furthermore, the type of~$A$ in the Cartan-Killing
nomenclature is uniquely determined by the cluster
algebra~$\AA$, and is called the ``cluster type'' of~$\AA$.
\end{theorem}

The study of cluster algebras of finite type reduces to the case
when the corresponding Cartan matrix~$A$ is indecomposable, i.e.,
is either of one of the classical types $A_n, B_n, C_n$, or~$D_n$,
or of one of the exceptional types $E_6, E_7, E_8, F_4$, or~$G_2$.
We refer to a seed with the exchange matrix~$B(A)$ as a
\emph{distinguished} seed of the corresponding cluster algebra.

\begin{theorem}[\cite{ca2}]
\label{th:conjectures-finite}
Parts~(1)-(3) of Conjecture~\ref{con:exchange-general} hold
for any cluster algebra~$\AA$ of finite type.
Furthermore, Conjectures~\ref{con:denominators}
and \ref{con:cluster-positive}~(1) hold for a distinguished
seed in~$\AA$.
\end{theorem}

The proof of Conjecture~\ref{con:denominators} (for a distinguished seed)
given in \cite[Theorem~1.9]{ca2} establishes a connection between
the cluster algebra structure and the corresponding root system.
Let $\Phi$ be the root system with the Cartan
matrix~$A$, and let $Q$ be the root lattice generated by~$\Phi$.
We identify~$Q$ with $\ZZ^n$ using the basis~$\Pi$ of simple
roots in~$\Phi$.
Let~$\Phi_{> 0}$ be the set of positive roots associated to~$\Pi$.
In this notation, \cite[Theorem~1.9]{ca2} combined with \cite[Theorem~1.8]{yga}
can be stated as follows.

\begin{theorem}
\label{th:denominator-roots}
The denominator vector parametrization
(see \eqref{eq:Laurent-normal-form}--\eqref{eq:denominator-vector})
with respect to a distinguished cluster~$\xx$ gives a bijection
between the set of cluster variables $\mathcal{X}$ and the set
$\Phi_{\geq -1} = \Phi_{> 0} \cup (- \Pi)$ of ``almost positive roots."
This bijection extends to a bijection between the set of all cluster
monomials and the root lattice~$Q = \ZZ^n$.
\end{theorem}

The first assertion of Theorem~\ref{th:denominator-roots} identifies
the clusters with certain subsets of~$\Phi_{\geq -1}$.
In view of \eqref{eq:delta-multiplicative}, the second one implies
that each of these subsets is a $\ZZ$-basis of the root lattice~$Q$,
and that the cones generated by them form a complete simplicial
fan in~$Q_\RR$.
As shown in \cite{cfz}, this fan is the normal fan
of a simple $n$-dimensional convex polytope
(the \emph{generalized associahedron} of
the corresponding type).
Thus, the exchange graph can be realized as the $1$-skeleton
of the generalized associahedron.

For the classical types, the clusters and generalized assohiahedra
have a nice concrete combinatorial realization in terms of regular
polygons and their triangulations; this realization is discussed
in detail in \cite{yga, cfz} and further explored in~\cite{ca4}.
For instance, it turns out that in type~$A_n$, the generalized
associahedron is the classical associahedron going back to
Stasheff's work~\cite{stasheff}, while in type~$B_n$ or~$C_n$,
this is the \emph{cyclohedron} appearing in \cite{bott-taubes}.
Using this realization we can prove the following.

\begin{theorem}[\cite{ca4}]
\label{th:conjectures-classical}
Conjectures~\ref{con:exchange-general},
\ref{con:denominators} and \ref{con:cluster-positive} hold
for any cluster algebra of classical type.
\end{theorem}

Having in mind Theorem~\ref{th:conjectures-classical}, it would be nice
to have a solution of the following tantalizing problem.

\begin{problem}
\label{prob:concrete-exceptional}
Find explicit combinatorial models for the generalized
associahedra of exceptional types.
\end{problem}

We conclude with a few remarks.

\begin{remark}
Given a matrix~$B$, it is a non-trivial problem to find a way to
decide whether the corresponding cluster algebra is of finite type.
This problem was solved by A.~Seven in \cite{seven-2}.
\end{remark}

\begin{remark}
Recall that our motivating example of an upper cluster algebra is
the coordinate ring $\CC[G^{u,v}]$ of a double Bruhat cell.
As shown in Section~\ref{sec:lecture1}, it is associated with a seed
$(F_\ii, \tilde B(\ii))$ given by \eqref{eq:Fi} and
\eqref{eq:tildeB-entries}.
The mutation procedure makes it obvious that, for every seed
$(\tilde \xx, \tilde B)$ mutation-equivalent to
$(F_\ii, \tilde B(\ii))$, the subset
$\tilde \xx \subset \CC[G^{u,v}]$ satisfies conditions (1) and (3)
in Definition~\ref{def:tp-basis}, and so provides global
coordinates in the totally positive variety $G^{u,v}_{> 0}$.
\end{remark}

\begin{remark}
It would be interesting to classify the double Bruhat cells such
that the corresponding cluster algebra is of finite type.
Some examples of this kind were given in \cite{ca3}.
An interesting feature of these examples is that
the cluster type of~$\CC[G^{u,v}]$ can be very different from
the Cartan-Killing type of the group~$G$.
For instance, as shown in \cite[Example~2.18]{ca3}, the cluster
algebra structure in $\CC[SL_3^{w_\circ,w_\circ}]$
(see Examples~\ref{example:SL3-w0w0} and \ref{ex:sl3-exchanges} above)
is of cluster type~$D_4$.
Another instance of the same phenomenon was demonstrated by
J.~Scott in \cite{scott}.
Namely, he introduced a natural cluster algebra
structure in the homogeneous coordinate ring of the
Grassmannian~${\rm Gr}_{k,n}$ of $k$-dimensional subspaces
of~$\CC^n$; among these cluster algebras, the following happen to be of finite
type: $\CC[{\rm Gr}_{2,n+3}]$ has cluster type~$A_n$ (as shown
in \cite{ca2}), while $\CC[{\rm Gr}_{3,6}]$ (resp.$\CC[{\rm
Gr}_{3,7}]$; $\CC[{\rm Gr}_{3,8}]$) has cluster type~$D_4$
(resp.~$E_6$; $E_8$).
\end{remark}

\vfill

\eject
%\pagebreak

\end{document}